\documentclass[12pt]{article}
\usepackage[english]{babel}
\usepackage{amssymb,amsmath,graphics}

\def\<{\langle}
\def\>{\rangle}

\newfont\recht{cmr12}

\def\cG{{\cal G}}

\def\cP{{\cal P}}

\def\cO{{\cal O}}
\def\cM{{\cal M}}

\def\qed{{\hfill\hphantom{.}\nobreak\hfill$\Box$}}
\def\proof{{\bf PROOF. }}
\newtheorem{lemma}{Lemma}[section]
\newtheorem{cor}[lemma]{Corollary}
\newtheorem{fact}[lemma]{Fact}
\newtheorem{remark}[lemma]{Remark}
\newtheorem{theo}[lemma]{Theorem}
\newtheorem{prop}[lemma]{Proposition}

\newtheorem{example}[lemma]{Example}
\newtheorem{question}[lemma]{Question}
\newtheorem{con}[lemma]{Conjecture}
\def\qed{{\hfill\hphantom{.}\nobreak\hfill$\Box$}}
\def\proof{{\bf Proof. }}
\newcommand{\ZZ}{\mathbb Z}
\parskip\medskipamount

\begin{document}
\title{\bf On finite and elementary generation of $SL_2(R)$}
\author{{\Large Peter Abramenko}}
\date{Department of Mathematics, P.O.Box
400137, University of Virginia, Charlottesville, VA 22904, USA, {\tt
pa8e@virginia.edu}}
\maketitle

\begin{abstract}
Motivated by a question of A. Rapinchuk concerning general reductive
groups, we are investigating the following question: Given a finitely
generated integral domain $R$ with field of fractions $F$, is there
a \emph{finitely generated subgroup} $\Gamma$ of $SL_2(F)$ containing
$SL_2(R)$? We shall show in this paper that the answer to this question is
negative for any polynomial ring $R$ of the form $R = R_0[s,t]$, where
$R_0$ is a finitely generated integral domain with infinitely many
(non--associate) prime elements. The proof applies Bass--Serre theory
and reduces to analyzing which elements of $SL_2(R)$ can be generated
by elementary matrices with entries in a given finitely generated
$R$--subalgbra of $F$. Using Bass--Serre theory, we can also exhibit new
classes of rings which do not have the $GE_2$ property introduced by P.M.
Cohn.
\end{abstract}

{\footnotesize {\bf Mathematics Subject Classification 2000:}
20H05, 20E06, 20E08; 13F20.  \\
{\bf Key words and phrases:} $SL_2$ over finitely generated domains,
$GE_2$, amalgams, Bass--Serre theory, Bruhat--Tits tree.}

\section{Introduction}

The starting point of the present paper was the following problem raised
by Andrei Rapinchuk (see \cite{Rap}, last paragraph):

\begin{question} \label{RP} \rm
Given a finitely generated integral domain $R$ with field of fractions $F$
and a reductive $F$--group $\cG$, does there exist a \emph{finitely generated}
subgroup of $\cG(F)$ which contains $\cG(R)$?
\end{question}

The background of this question is the following. In \cite{RSS} Rapinchuk,
Segev and Seitz prove the beautiful theorem that any finite quotient of
the multiplicative group of a finite dimensional division algebra $D$ is
solvable. This leads them to the question whether any finite quotient of
$\cG(F)$ for a reductive group $\cG$ over an infinite field $F$ is solvable.
Their result shows that this is true for $\cG = GL_{1,D}$, and the obvious
next candidate is $\cG = SL_{1,D}$. However, the transition from
$GL_{1,D}$ to $SL_{1,D}$ is non--trivial and involves the question whether
a normal subgroup of finite index in $SL_{1,D}(F)$ contains a finite index
subgroup which is normal in $GL_{1,D}(F)$. Rapinchuk could answer this
question in the affirmative provided that Question~\ref{RP} has a positive
answer for $\cG = SL_{1,D}$.

However, even without this background, Question~\ref{RP} is interesting and
challenging. It is certainly well--known that it has a positive answer for
$S$-arithmetic groups $\cG(R)$, which in "almost all" cases are finitely
generated themselves (see Remark~\ref{Nagao} below for more precise statements).
But it was not clear whether one could expect a positive answer to
Question~\ref{RP} for arbitrary finitely generated integral domains $R$.
A standard reduction in algebraic $K$-theory provides, modulo a (difficult)
problem concerning the finite generation of $K_1(R)$ for regular $R$, some
evidence that Rapinchuk's question admits a positive answer for
$\cG(R) = SL_n(R)$ in case $n$ is "sufficiently large" (see Question~\ref
{K1} and Remark~\ref{n-big}).

On the other hand, it seemed unlikely to me that Question~\ref{RP} had a
positive answer for $\cG = SL_2$. The present paper is (mainly) about turning
this vague idea into a rigorous proof for a reasonable class of rings $R$.
This is the following theorem which will be proved in Section 4 (see
Theorem~\ref{Mainresult}).

\begin{theo} \label{mainresult}
Let $R_0$ be a finitely generated integral domain with infinitely many
non--associate primes, $R = R_0[s,t]$ with field of fractions $F$ and
$\Gamma$ a group with $SL_2(R) \leq \Gamma \leq SL_2(F)$. Then
$\Gamma$ is not finitely generated.
\end{theo}

The strategy is the following. One starts by making the elementary observation
that a finitely generated group $\Gamma$ with $SL_2(R) \leq \Gamma \leq SL_2(F)$
exists if and only if $SL_2(R) \subseteq E_2(S)$ for some finitely generated
$R$--subalgebra $S$ of $F$ (see Lemma~\ref{reformulation}; $E_2(S)$ denotes the
subgroup of $SL_2(S)$ generated by elementary matrices). So for any given $S$,
one wants to exhibit an element of $SL_2(R)$ which is not contained in $E_2(S)$.
To this end, one provides $F$ with an appropriate valuation, let $G = SL_2(F)$
act on the corresponding (Bruhat--Tits) tree $T$, and considers the subgroup
$H_0$ of $H =SL_2(S)$ generated by the stabilizers in $H$ of the two vertices
of a fundamental edge of $T$. By another elementary observation (Lemma~
\ref{1/r}), $H_0$ contains $E_2(S)$. Now Bass-Serre theory provides us with
criteria to decide whether $H = H_0$ and with a method to construct a
\emph{concrete} element $h$ of $H$ not in $H_0$ if $H \neq H_0$ (see
Lemmas~\ref{L1} -- \ref{L3}). These general criteria now have to be applied in
the given situation (which also requires a bit of commutative algebra), finally
yielding matrices $h \in SL_2(R)$ with $h \not\in E_2(S)$.

It turns out that this method is also effective in order to establish, under
certain conditions, that $SL_2(R[1/\pi]) \neq E_2(R[1/\pi])$ for a (not
necessarily finitely generated) integral domain $R$ with prime element $\pi$.
More precisely, we obtain the following theorem which will be proved, among
other results, in Section 3 (see Corollary~\ref{notGE2}).

\begin{theo} \label{NotGE2}
Let $R$ be an integral domain and $\pi$ a prime element of $R$ satisfying
$\bigcap_{n \geq 0}\pi^nR$ $ = \{0\}$. If $R/\pi R$ is not a Bezout domain or
if the canonical homomorphism $SL_2(R) \rightarrow SL_2(R/\pi R)$ is not
surjective, then $SL_2(R[1/\pi]) \neq E_2(R[1/\pi])$.
\end{theo}

This generalizes results about Laurent polynomial
rings proved in \cite{BM} and \cite{Chu}. Here (that is in
Theorem~\ref{Equivalent}) we again investigate the
question whether $H = H_0$, with $H = SL_2(R[1/\pi])$ acting on the
Bruhat--Tits tree associated to $SL_2(F)$ and $\pi$, and $E_2(R[1/\pi]) \leq
H_0$. In fact we can derive a necessary and sufficient condition for $H = H_0$
in this situation.

The paper is organized as follows. In Section 2 we prove some lemmas about
subgroups of amalgams, using the action of these groups on the associated trees.
This provides us with the above mentioned criteria concerning $H = H_0$,
$H \neq H_0$. We first apply these criteria in Section 3 in order to deduce a
necessary and sufficient condition for $SL_2(R[1/\pi]) = H_0$, where $H_0$ is
the subgroup generated by $SL_2(R)$ and its conjugate by the diagonal matrix
with entries $1/\pi$ and $1$. In Section 4 we deduce the negative answer to
Rapinchuk's problem for the groups $SL_2(R_0[s,t])$ in the way indicated
above. We conclude this paper by listing some further questions and
conjectures in Section 5.

{\bf Acknowledgements:} I am grateful to my colleagues Andrei Rapinchuk and
Nicholas Kuhn for stimulating discussions on this subject. Andrei's problem for
reductive groups motivated the whole paper, and Section~\ref{notGE-2} emerged
from an interesting question which Nick asked me during one of my talks.

\section{About subgroups of amalgams}\label{amalgams}

In this section, we consider the following set--up. The group $G$ is the free
product with amalgamation of its two subgroups $A$ and $B$, amalgamated along
their intersection $U = A \cap B$, $H$ is an arbitrary subgroup of $G$ and $H_0$
is the subgroup of $H$ generated by $(A \cap H) \cup (B \cap H)$.
We use the notations
$$G = A \, *_U \, B \quad , \quad H \leq G \quad {\rm and} \quad
H_0 = \langle (A \cap H) \cup (B \cap H) \rangle$$
We are interested in the following

\begin{question}\label{qu} \rm
When is $H = H_0$, and when is $H \neq H_0$?
\end{question}

Structure theorems for subgroups of amalgams have been known in combinatorial
group theory for a long time, see for instance \cite{KS}. However,
Question~\ref{qu} is attacked in a less technical and more transparent way
by using group actions on trees. If $X$ is
a tree, we denote by $VX$ its set of vertices and by $EX$ its set of edges.
Here an edge is always understood as a geometric edge, i.e. it is identified
with a subset of cardinality 2 of $VX$. If a group $C$ acts (on the left) on
$X$, then we denote this action with a dot and set
$C_{\alpha} := \{ c \in C \mid c.\alpha = \alpha \}$ for any vertex or any
edge $\alpha$ of $X$.
Let us first recall one of the basic results about amalgams
(cf. \cite[Chapter I, Section 4.1]{Se}).

\begin{fact}\label{fact1}
$G$ acts without inversion on a (suitable) tree $T$ with an edge $e = \{x,y\}$ as
fundamental domain (i.e. $G.e = ET$, $G.x \cup G.y = VT$ and $G.x \neq G.y$)
such that $G_x = A$, $G_y = B$ and $G_e = U$.\\
Conversely, if a group $G'$ acts without inversion on a tree $T'$ with an edge
$e' = \{x',y'\} \in ET'$ as fundamental domain, then
$G' = G_{x'} \, *_{G'_{e'}} \, G'_{y'}$.
\end{fact}

The second statement in Fact~\ref{fact1}  has a well-known generalization
(cf. \cite[Chapter I, Section 4.5]{Se}).

\begin{fact}\label{fact2}
If a group $G'$ acts on a tree $T'$ with a subtree $T_1$ of $T$ as fundamental
domain, and if $T_1$ is considered as a tree of groups with respect to the
system ${\cal G} = \bigl( (G'_v)_{v \in VT_1}$, $(G'_f)_{f \in ET_1} \bigr)$,
then $G'$ is canonically isomorphic to the direct limit (which is an "amalgam
along $T_1$") ${\rm lim}({\cal G}, T_1)$ of this tree of groups.
\end{fact}

We now fix $T$ and $e$ as in Fact~\ref{fact1}.
In the following sequence of three lemmas dealing with Question~\ref{qu},
the first one is similar to some well--known results. However, for the
convenience of the reader I shall give a short proof also in this case.

\begin{lemma}\label{L1}
Denote by $X$ the subforest of $T$ with edge set $H.e$ and vertex set
$H.x \cup H.y$.
\begin{itemize}
\item[(i)] If $X$ is connected, then
$H = (A \cap H) *_{(U \cap H)} (B \cap H)$, and in particular $H = H_0$.
\item[(ii)] If $A = (A \cap H)U$ and $B = (B \cap H)U$, then $X = T$,
and hence $H = H_0$ by (i).
\end{itemize}
\end{lemma}

\proof Observe that $H$ acts on $T$ with stabilizers
$H_x = G_x \cap H = A \cap H$, $H_y = B \cap H$ and $H_e = U \cap H$.
Now (i) immediately follows from the second part of Fact~\ref{fact1} if we
set $G' = H$ and $T' = X$, which is a tree by assumption.

In order to prove (ii), we have to show $H.e = ET$. Given $e' \in ET$,
we consider the geodesic
$\gamma = (z = z_0,z_1,\ldots,z_n)$ in $T$ with $z_i \in VT$ for all $i$,
$\{z_i,z_{i+1}\} \in ET$ for all $i \leq n-1$, $z_i \neq z_{i+2}$ for all
$i \leq n-2$, $\{z_0,z_1\} = e$ and
$\{z_{n-1},z_n\} = e'$. We show $e' \in H.e$ by induction on $n$. We may
assume $z_1 = x$ (the case $z_1 = y$ is similar) and $n \geq 2$. By
assumption, $H_xG_e = G_x$. Now $G_x$ acts transitively on the set of edges
containing $x$ because $e$ is a fundamental domain for the action of
$G$ on $T$. Hence there exists an $h \in H_x$ such that $h.\{z_1,z_2\} = e$.
Applying the induction hypothesis to the geodesic
$(x = h.z_1,y = h.z_2,\ldots,h.z_n)$, we obtain $h.e' \in H.e$, which
immediately implies $e' \in H.e$.
\qed

If one of the two assumptions in Lemma~\ref{L1}(ii) is not satisfied, then $H$
is "often" different from $H_0$, as the following result shows.

\begin{lemma}\label{L2}
Assume that the following two conditions are satisfied.
\begin{itemize}
\item[(1)] There exists $a \in A$ with $a \not\in (A \cap H)U$.
\item[(2)] There exists $b \in B$ with $b \not\in U$ and $aba^{-1} \in H$.
\end{itemize}
Then $H \neq H_0$.
\end{lemma}

\proof Note first that the edges $e$ and $a.e$ are in different $H$--orbits
since $a.e \in H.e$ implies $a \in HU \cap A = (A \cap H)U$, contradicting (1).
Because the element $h := aba^{-1} \in aBa^{-1} = G_{a.y}$ is also in $H$
but not in $aUa^{-1} = G_{a.e}$ by assumption (2), we have
$h \in H_{a.y}$ and $h \not\in H_{a.e}$. We now distinguish two cases.

{\bf First Case: $a \in HB$.} Here $a.y$ is contained in $H.y$. So the images of
$y$ and $a.y$ are equal in the quotient graph $H \setminus T$. Since $a \in G_x$,
we also have $x = a.x$. But as observed above, the edges $e$ and $a.e$ do not
have the same image in $H \setminus T$. Therefore $H \setminus T$ contains a
circuit of length 2. It now follows from \cite[Chapter I, Section 5.4,
Corollary 1 of Theorem 13]{Se} that $H \neq \langle
\bigcup_{v \in VT}H_v \rangle$, hence in particular $H \neq H_0$.

{\bf Second Case: $a \not\in HB$.} So the $H$--orbits of $a.y$ and $y$ are
different. Hence the subtree $T_0$ of $T$ with $VT_0 = \{y,x,a.y\}$ and
$ET_0 = \{e,a.e\}$ is mapped injectively into $H \setminus T$ by the
canonical projection $T \rightarrow H \setminus T$. Now consider the
subgroup $H' := \langle H_y \cup H_x \cup H_{a.y} \rangle$ of $H$ which
acts on $T' := H'.T_0$. We show that $T'$ is connected and hence a subtree
of $T$. For any integer $n \geq 1$, we set
$T_n := \bigcup h_1\ldots h_n.T_0$, where $h_1\ldots h_n$ runs over all
products with factors $h_i \in (H_y \cup H_x \cup H_{a.y})$ for all
$1 \leq i \leq n$. For any such product, the intersection
$h_1\ldots h_n.T_0 \cap T_{n-1}$ is obviously nonempty. So by induction,
$T_n$ is connected for all $n$. Hence also $T' = \bigcup_{n\geq 0}T_n$
is connected. By construction and since $T_0$ embeds into $H \setminus T$, hence
also into $H' \setminus T'$, $T_0$ is a fundamental domain for the action of
$H'$ on $T'$. So by Fact~\ref{fact2}, $H'$ is the direct limit of the tree of
groups associated with $T_0$ and $\bigl( (H_y,H_x,H_{a.y}), \, (H_e,H_{a.e})
\bigr)$, showing
$$H' = (H_y \, *_{H_e} \, H_x) \, *_{H_{a.e}} \, H_{a.y}
     =  H_0 \, *_{H_{a.e}} \, H_{a.y}$$
As observed above, $h \in H_{a.y}$ and $h \not\in H_{a.e}$. The normal form for
amalgams (cf. \cite[Chapter I, Section 1.2]{Se})  now yields $h \not\in H_0$.
Therefore $H \neq H_0$.
\qed

The proof of Lemma~\ref{L2} yields additional information which is worth mentioning.

\begin{lemma}\label{L3}
Assume that the following two conditions are satisfied.
\begin{itemize}
\item[(1)] There exists $a \in A$ with $a \not\in HB$.
\item[(2)] There exists $b \in B$ with $b \not\in U$ and $aba^{-1} \in H$.
\end{itemize}
Then $h := aba^{-1}$ is not an element of $H_0$.
\end{lemma}

\proof $a \not\in HB$ obviously implies $a \not\in (A \cap H)U$.
So the assumptions of Lemma~\ref{L2} are satisfied, and additionally we are in
Case 2 of its proof. As demonstrated there, this implies $h \not\in H_0$.
\qed.

In the application which we shall discuss in the last section it will become
important that Lemma~\ref{L3} provides us with a method that produces
\emph{concrete} elements in $H$ which are not contained in $H_0$.
It turns out that Condition (2) can be trivially satisfied in those situations
where we are going to apply Lemma~\ref{L2} and Lemma~\ref{L3}. However, some
work will be necessary in order to verify Condition (1).

\section{Non--elementary generation of $SL_2(R[1/{\pi}])$}
\label{notGE-2}

Let $R$ be an integral domain with field of fractions $F$
and $\pi \in R$ a prime element.
In this section we shall deduce some necessary conditions for $SL_2(R[1/{\pi}])$
to be generated by elementary matrices. We start with an easy exercise in
commutative algebra which we shall need later on.

\begin{lemma}\label{lprincipal}
Let $u,v,x,y \in R$ with $yu \neq 0$ and $ux = vy$. Then $(u,v)$ is a principal
ideal of $R$ if and only $(x,y)$ is.
\end{lemma}

\proof By symmetry we may assume that $(u,v) = (d)$ with $d \in R$.
Then $d \neq 0$ (since $u \neq 0$), and $u_1 := u/d \in R$ as well as
$v_1 := v/d \in R$. Also there exist $r,s \in R$ with $ru + sv = d$,
hence $ru_1 + sv_1 = 1$. We claim that $(x,y) = (sx + ry)$. So we have
to show $x,y \in (sx + ry)$. Recall that $ux = vy$, hence $u_1x = v_1y$.
So we obtain $x = (ru_1 + sv_1)x = rv_1y + sv_1x = v_1(sx + ry)$ and
$y = (ru_1 + sv_1)y = ru_1y + su_1x = u_1(sx + ry)$. This proves the claim.
\qed



With respect to elementary matrices, we shall use the following notations. We set
$$E_{12}(r) := \left( \begin{array}{rr} 1 & r\\ 0 & 1\\ \end{array} \right) , \quad
E_{21}(r) := \left( \begin{array}{rr} 1 & 0\\ r & 1\\ \end{array} \right)  \quad
{\rm for} \ {\rm any} \  r \in R$$
and then $E_{12}(R) := \{ E_{12}(r) \mid r \in R \}$,
$E_{21}(R) := \{ E_{21}(r) \mid r \in R \}$ as well as
$E_2(R) := \langle E_{12}(R) \cup E_{21}(R) \rangle $. For any two $\alpha, \beta
\in F^*$, we define
$D(\alpha,\beta) := \left( \begin{array}{rr} \alpha & 0\\ 0 & \beta\\ \end{array}
\right) \in GL_2(F)$.

The following observation concerning the ring $R[1/r]$ is elementary but useful.
It was already successfully applied in \cite{BM}.

\begin{lemma}\label{1/r}
For any $r \in R$ with $r \neq 0$ we obtain the following:
\begin{itemize}
\item[(i)] $E_2(R[1/r]) = \langle E_2(r) \cup \{ E_{12}(1/r) \} \rangle$
\item[(ii)] $E_2(R[1/r]) \leq \langle SL_2(R) \cup D(1/r,1)SL_2(R)D(r,1)
\rangle \leq SL_2(R[1/r])$
\end{itemize}
\end{lemma}

\proof (i) follows from the well--known identity
$$\large{(}E_{12}(1/r)E_{21}(-r)E_{12}(1/r)\large{)}\large{(}E_{12}(-1)
E_{21}(1)E_{12}(-1)\large{)} = D(1/r,r)$$ together with
$D(1/r,r)^nE_{12}(R)D(1/r,r)^{-n} = E_{12}(r^{-2n}R)$ and
$D(1/r,r)^{-n}E_{21}(R)D(1/r,r)^n$  $= E_{21}(r^{-2n}R)$ for all integers $n$.
The first inclusion in (ii) follows immediately from (i), and the second is
obvious.
\qed

So $SL_2(R[1/r])$ can only be generated by elementary matrices if it is also
generated by $SL_2(R) \cup D(1/r,1)SL_2(R)D(r,1)$. It is the question whether
$SL_2(R[1/r]) = \langle SL_2(R) \cup D(1/r,1)SL_2(R)D(r,1) \rangle$ or not to
which we can apply our results from Section~\ref{amalgams}. More precisely,
we shall do this in case $r = \pi$ is a prime element in $R$ in order to have
a nice action of $SL_2(R[1/r])$ on a suitable Bruhat--Tits tree (see
Fact~\ref{BTtree} below).
Before we can introduce the latter, we need the following assumption,
which is obviously satisfied for all noetherian, and hence also for all
finitely generated rings.

{\bf Assumption (A):} The prime element $\pi \in R$ satisfies
$\bigcap_{n \geq 0}\pi^nR = \{0\}$.


Now let (A) be satisfied for a fixed prime $\pi \in R$. We define a $\pi$--adic
valuation $v = v_{\pi}$ on $F$ in the usual way. For any $r \in R \setminus
\{0\}$, we set $v(r) := \max\{n \geq 0 \mid r \in (\pi^n)\}$, which exists
in view of (A). We further define $v(x/y) := v(x) - v(y)$ for $x,y \in
R \setminus \{0\}$ and $v(0) := \infty$. It is immediately verified that $v$
is thus a \emph{discrete valuation} on $F$. We denote by $\cO$ the associated
discrete valuation ring $\cO = \{\alpha \in F \mid v(\alpha) \geq 0 \}$, by
$\cP = \pi\cO$ its maximal ideal and by ${\cO}^* = \cO \setminus \cP$ its
group of units.
Our reference for the following statements is again Serre's book;
cf. \cite[Chapter II, Section 1]{Se}.

\begin{fact}\label{BTtree}
Given $F$ together with the disrete valuation $v$, one can construct a tree $T$
(which is also the Bruhat-Tits building of $SL_2(F)$ with respect to $v$)
on which $G = SL_2(F)$ acts without inversion and with an edge as fundamental
domain. This edge $e$ can be chosen such that the stabilizers of its two
vertices are $A = SL_2(\cO)$ and $B = D(1/\pi,1)SL_2(\cO)D(\pi,1)$, respectively,
and $G_e = U = A \cap B = SL_2\left( \begin{array}{rr} \cO & \cO\\ \cP &
\cO\\ \end{array} \right)$. By Fact~\ref{fact1} this implies
$G = A \, *_U \, B$.
\end{fact}

We now want to apply the results of Section~\ref{amalgams} to this situation.
Recall that a commutative ring is called Bezout if each of its finitely
generated ideals is a principal ideal.

\begin{theo}\label{Equivalent}
If Assumption (A) is satisfied, then the following two statements are equivalent:
\begin{itemize}
\item[(1)] $SL_2(R[1/\pi]) = \langle SL_2(R) \cup D(1/\pi,1)SL_2(R)D(\pi,1) \rangle$
\item[(2)] $R/\pi R$ is a Bezout domain and the canonical homomorphism
$SL_2(R) \rightarrow SL_2(R/\pi R)$ is surjective
\end{itemize}
\end{theo}

\proof We consider the subgroup $H = SL_2(R[1/\pi])$ of $G = SL_2(F)$. With
$A,B$ and $U$ as in Fact~\ref{BTtree}, we obtain $A \cap H = SL_2(R)$,
$B \cap H = D(1/\pi,1)(A \cap H)D(\pi,1) = D(1/\pi,1)SL_2(R)D(\pi,1)$ and
$U \cap H = SL_2\left( \begin{array}{rr} R & R\\ \pi R & R\\ \end{array} \right)$.
Setting $H_0 = \langle (A \cap H) \cup (B \cap H) \rangle =
\langle SL_2(R) \cup D(1/\pi,1)SL_2(R)D(\pi,1) \rangle$, we have to answer
Question~\ref{qu} in this situation.

\emph{ The implication "(2) $\Rightarrow$ (1)".} We shall show that
$A = (A \cap H)U$ and $B = (B \cap H)U$ if (2) is satisfied. Then (1) will
follow from Lemma~\ref{L1}. So, firstly, given any
$a = \left( \begin{array}{rr} \alpha & *\\ \beta & *\\ \end{array} \right)
\in A = SL_2(\cO)$, we have to find an
$h = \left( \begin{array}{rr} \ * & *\\ r & s\\ \end{array} \right) \in
A \cap H = SL_2(R)$ such that $ha \in U$, i.e. such that
$r\alpha + s\beta \in \cP$. Since $\alpha, \beta \in \cO$, there exists
$z \in R \cap {\cO}^*$ such that $p := z\alpha, q := z\beta$ are both elements
of $R$. Denote by ${\overline p}$ and ${\overline q}$ the respective images in
$R/\pi R$. Because this ring is Bezout by assumption, $({\overline p},{\overline q})$
is a prinicipal ideal, i.e. $({\overline p},{\overline q}) = (\delta)$ for some
$\delta \in R/\pi R$. Note that $\delta \neq {\overline 0}$ since
$\cO = \alpha \cO + \beta \cO = p\cO + q\cO$.
Set $\lambda := {\overline q}/\delta, \,  \mu := -{\overline p}/\delta \in R/\pi R$.
Then $(\lambda,\mu) = ({\overline 1})$, and we can find a matrix
$ \left( \begin{array}{rr} \ * & *\\ \lambda & \mu \\ \end{array} \right) \in
SL_2(R/\pi R)$. By assumption, this matrix has a preimage
$ \left( \begin{array}{rr} \ * & *\\ r & s\\ \end{array} \right)$ in $SL_2(R)$;
we call this preimage $h$. Now by construction,
${\overline r}{\overline p} + {\overline s}{\overline q} =
\lambda {\overline p} + \mu {\overline q} =
({\overline q}{\overline p} - {\overline p}{\overline q})/\delta = {\overline 0}$.
Therefore, $rp + sq \in \pi R$, and hence also $r\alpha + s\beta \in \cP$, because
$z \in {\cO}^*$. This proves that $ha \in U$.

The equation $B = (B \cap H)U$ is equivalent to
$A = (A \cap H)D(\pi,1)UD(1/\pi,1)$. Now this equation can be proved completely
similar as the equation $A = (A \cap H)U$ above. We only have to produce a
$(1,2)$--entry in $\cP$ for the product $ha$ instead of a $(2,1)$--entry in $\cP$.

\emph{ The implication "(1) $\Rightarrow$ (2)".} Now we assume that $R/\pi R$ is
not Bezout or that the canonical homorphism $\phi: SL_2(R) \rightarrow
SL_2(R/\pi R)$ is not surjective. In the first case, we choose elements
$x,y \in R$ such that the ideal $({\overline x},{\overline y})$ of $R/\pi R$ is
not principal. In the second case, we choose $x,y \in R$ such that there exists
a matrix $k \in SL_2(R/\pi R)$ of the form
$k = \left( \begin{array}{rr}
{\overline y} & *\\ {\overline x} & *\\ \end{array} \right)$ with
$k \not\in {\rm im}\, \phi$. Note that $x,y \in {\cO}^*$ in both cases. This is
obvious if $({\overline x},{\overline y})$ is not principal, and in the second case
it follows from $E_2(R/\pi R) \subseteq {\rm im} \, \phi$ and the easy observation
that $k \in E_2(R/\pi R)$ if one of the entries of $k$ is equal to ${\overline 0}$.

Having chosen $x$ and $y$, we now define
$a := \left( \begin{array}{rr} \ 1 & 0\\ x/y & 1 \\ \end{array} \right) \in A$.
Suppose $a \in (A \cap H)U$. Then, as above, there is an
$h = \left( \begin{array}{rr} \ * & *\\ r & s \\ \end{array} \right) \in SL_2(R)$
such that $ha \in U$. This implies $r + sx/y \in \cP$, hence
$ry + sx \in \cP \cap R = \pi R$ and thus
${\overline r}{\overline y} + {\overline s}{\overline x} = {\overline 0}$ in
$R/\pi R$. However, this contradicts our choice of $x$ and $y$ in both cases.
Firstly, ${\overline r}{\overline y} = - {\overline s}{\overline x}$ and
Lemma~\ref{lprincipal} imply that $({\overline x}, {\overline y})$ is principal
because $({\overline r},-{\overline s}) = ({\overline 1})$ is principal.
Secondly, $\phi(h)k$ is of the form $\phi(h)k =
\left( \begin{array}{rr} \ * & *\\ {\overline 0} & * \\ \end{array} \right)$.
Hence $\phi(h)k \in E_2(R/\pi R) \subseteq {\rm im} \, \phi$, implying
$k \in {\rm im} \, \phi$.

So in both cases, $a \in A$ and $a \not\in (A \cap H)U$. We now set
$b := \left( \begin{array}{rr} \ 1 & y^2/\pi \\ 0  & 1 \\ \end{array} \right)
\in B$. We have $b \not\in U$ since $y \in {\cO}^*$. And we have that
$aba^{-1} = \left( \begin{array}{rr}
1 - xy/\pi & y^2/\pi \\ -x^2/\pi & 1 + xy/\pi\\ \end{array} \right)
\in SL_2(R[1/\pi] = H$. Therefore, $H \neq H_0$ by Lemma~\ref{L2}.
\qed
\vskip5mm

From my point of view, the most interesting consequence of
Theorem~\ref{Equivalent} (and of the elementary Lemma~\ref{1/r})
is the following, which was stated as Theorem~\ref{NotGE2} in the
Introduction.

\begin{cor}\label{notGE2}
If $R$ is an integral domain and $\pi \in R$ a prime element satisfying
Assumption (A), then $SL_2(R[1/\pi]) \neq E_2(R[1/\pi])$ whenever $R/\pi R$
is not Bezout or the canonical homomorphism $SL_2(R) \rightarrow
SL_2(R/\pi R)$ is not surjective.
\end{cor}
\qed

\begin{remark}\label{noetherian} \rm
If $R$ is noetherian, then Assumption (A) is automatically satisfied and
hence superfluous in the statement of Theorem~\ref{Equivalent} as well as
in Corollary~\ref{notGE2}. Furthermore,  "Bezout" can be equivalently
replaced with "principal ideal domain" in this case.
\end{remark}

A special case of Corollary~\ref{notGE2} is obtained if $R = R_0[t]$ is the
polynomial ring in one variable over an integral domain $R_0$ and $\pi = t$,
in which case Assumption (A) is clearly satisfied and the canonical
homorphism $SL_2(R) \rightarrow SL_2(R/\pi R) = SL_2(R_0)$ always surjective.
So we recover the following result about Laurent polynomial rings which was
partly deduced by Bachmuth-Mochizuki in \cite {BM}
and first proved in the generality we state it here by H. Chu (cf. \cite{Chu}).

\begin{cor}\label{Laurent}
If $R_0$ is an integral domain which is not Bezout, then
$SL_2(R_0[t,t^{-1}]) \neq E_2(R_0[t,t^{-1}])$.
\end{cor}
\qed

One remarkable feature about Corollary~\ref{notGE2} is that $R[1/\pi]$ cannot
be a $GE_2$--ring in the sense of Cohn (cf. \cite{Co}), no matter how "nice"
$R$ is, if the quotient ring $R/\pi R$ does not have the stated properties.
So, in particular, the $GE_2$--property is not preserved by the process of
"localization", by which we mean the transition from a ring to one of its
rings of fractions. This is well demonstrated by the following example and
answers a respective question of Nick Kuhn.

\begin{example}\label{local} \rm
Let $R$ be a (noetherian) regular local ring of Krull dimension $\geq 3$.
(Take for instance the localization $R = S_{\cM}$ of the polynomial ring
$S = K[t_1,t_2,t_3]$ over a field $K$ at the maximal ideal
$\cM = (t_1,t_2,t_3)$.) Because $R$ is local, $SL_2(R) = E_2(R)$.
However, for any prime element $\pi \in R$ (and $R$ has a
lot of prime elements since it is a unique factorization domain),
$R/\pi R$ has Krull dimension $\geq 2$. Hence $R/\pi R$ is not a principal
ideal domain, and so $SL_2(R[1/\pi]) \neq E_2(R[1/\pi])$ by
Corollary~\ref{notGE2}.
\end{example}

So far we have been discussing consequences of Theorem~\ref{Equivalent}
concerning the elementary generation of $SL_2(R[1/\pi])$. Let us finish
this section by mentioning two cases where Condition (2) is obviously
satisfied.

\begin{cor}\label{Dedekind}
Let $R$ be a ring which is either a Dedkind domain or of the form
$R = R_0[t]$ with a Bezout domain $R_0$. Let $\pi$ be an arbitrary
prime element of $R$ in the first case and $\pi = t$ in the second case.
Then $SL_2(R[1/\pi]) = \langle SL_2(R) \cup D(1/\pi,1)SL_2(R)D(\pi,1)
\rangle$, and moreover
$$SL_2(R[1/\pi]) = SL_2(R) \, *_U \,
SL_2\left( \begin{array}{rr} \ R & \pi^{-1}R \\ \pi R & R \\ \end{array} \right)
 \ {\rm with} \
U = SL_2\left( \begin{array}{rr} \ R & R\\ \pi R & R \\ \end{array} \right)$$
\end{cor}

\proof $R/\pi R$ is a field in the first, and $R/\pi R = R_0$ in the second
case. Hence Condition (2) of Theorem~\ref{Equivalent} is satisfied, yielding
the first claim of this corollary. However, the proof of Theorem~\ref{Equivalent}
in fact shows that the assumptions of Lemma~\ref{L1}(ii) are satisfied.
Therefore, this lemma implies the second
claim about the amalgam presentation of $SL_2(R[1/\pi])$.
\qed

For Dedekind rings, a different proof of Corollary~\ref{Dedekind} is given
in \cite[Chapter II, Section 1.4]{Se}. (It is stated there only for
$R = {\mathbb Z}$ but could be generalized.) The result about Laurent
polynomial rings is essentially Theorem 2 in \cite{BM}.

\begin{remark} \rm
Whenever Condition (2) of Theorem~\ref{Equivalent} is satisfied, the proof
of the implication "(2) $\Rightarrow$ (1)" together with Lemma~\ref{L1}
yields the same amalgam presentation of $SL_2(R[1/\pi])$ as stated in
Corollary~\ref{Dedekind}.
\end{remark}

\begin{remark} \label{ZLaurent} \rm
Corollary~\ref{Dedekind} is not of much help in order to decide the question
whether Laurent polynomial rings in one variable over principal ideal domains
are elementary generated. To the best of my knowledge, it is still an open
problem whether the groups $SL_2({\mathbb Z}[t,t^{-1}])$ and
$SL_2({\mathbb F}_q[t_1,t_1^{-1};t_2,t_2^{-1}])$ are generated by elementary
matrices. Even the weaker question whether they are finitely generated still
seems to be open. (It is easily seen that $E_2({\mathbb Z}[t,t^{-1}])$ and
$E_2({\mathbb F}_q[t_1,t_1^{-1};t_2,t_2^{-1}])$ are finitely generated;
see the proof of Lemma~\ref{elementary}(ii) below.)
\end{remark}

\section{Non--finite generation of groups between $SL_2(R)$ and $SL_2(F)$}
\label{infinite}

We now turn to the problem which motivated this paper. In this section, $R$
will always denote a finitely generated integral domain, i.e. an integral domain
which is finitely generated as a ${\mathbb Z}$--algebra. So $R$ can be obtained
by adjoining finitely many elements to its prime ring $P$
($P = {\mathbb Z}$ or $P = {\mathbb F}_p$), that is
$R = P[x_1,\ldots,x_n]$ with elements $x_1,\ldots,x_n \in R$. We denote by $F$
the field of fractions of $R$, so $F = Q(x_1,\ldots,x_n)$ with
$Q = {\mathbb Q}$ or $Q = {\mathbb F}_p$. Let us recall the question we
want to answer:

\begin{question}\label{Qu} \rm
Does there exist a \emph{finitely generated} group $\Gamma$ with
$SL_2(R) \leq \Gamma \leq SL_2(F)$?
\end{question}

There is an intimate connection between finite and elementary generation,
as the following easy lemma shows.

\begin{lemma} \label{elementary} The following holds:
\begin{itemize}
\item[(i)] Any finitely generated subgroup of $SL_2(F)$ is contained in
$E_2(S)$ for some finitely generated subring $S$ of $F$.
\item[(ii)] Any finitely generated subring $S$ of $F$ is included in
some finitely generated $S$--subalgebra $S' \subseteq F$ for which
$E_2(S')$ is a finitely generated group.
\end{itemize}
\end{lemma}

\proof We again denote by $P = {\mathbb Z}$, respectively
$P = {\mathbb F}_p$, the prime subring of $F$.

(i) Assume that $\Gamma = \langle \gamma_1,\ldots,\gamma_l \rangle$
is a finitely generated subgroup of $SL_2(F)$. Recall that
$SL_2(F) = E_2(F)$. For each $1 \leq i \leq l$, we fix a representation
of $\gamma_i$ as a product $\gamma_i = e_{i1}\ldots e_{ik_i}$ of
elementary matrices $e_{ij}$ ($\ 1 \leq j \leq k_i$) in $SL_2(F)$. Let
$M_i$ be the finite subset of $F$ consisting of all entries of all the
$e_{ij}$. Set $M := \bigcup_{i\leq l}M_i$ and $S := P[M]$. Then, by
construction, $S$ is a finitely generated subring of $F$ and
$\Gamma = \langle \gamma_1,\ldots,\gamma_l \rangle \leq E_2(S)$.

(ii) Now suppose that $S = P[y_1,\ldots,y_m]$ with $y_i \in F^*$ for all
$1 \leq i \leq m$. We set $S' := S[y_1^{-1},\ldots,y_m^{-1}]$. An
easy calculation (using conjugation of elementary matrices by diagonal
matrices) shows that $E_2(S')$ is generated by the diagonal matrices
$D(y_i,y_i^{-1})$, $1 \leq i \leq m$, together with the elementary matrices
$E_{12}(z)$, $E_{21}(z)$, where $z$ runs over all products of the form
$z = y_{i_1}\ldots y_{i_k}$ with $1 \leq i_1 < i_2 < \ldots < i_k \leq m$,
including the empty product ($k = 0$) which is 1 by definition. In
particular, $E_2(S')$ is finitely generated.
\qed

This admits a reformulation of Question~\ref{Qu} in terms of elementary
generation.

\begin{lemma} \label{reformulation}
With $R$ and $F$ as above, the following two statements are equivalent:
\begin{itemize}
\item[(1)] There exists a finitely generated group $\Gamma$ with
$SL_2(R) \leq \Gamma \leq SL_2(F)$.
\item[(2)] There exists a finitely generated $R$--subalgebra
$S \subseteq F$ such that $SL_2(R) \leq E_2(S)$.
\end{itemize}
\end{lemma}

\proof If (1) is satisfied, then $\Gamma$ and hence $SL_2(R)$ is contained
in $E_2(S)$ for a finitely generated subring $S$ of $F$ by
Lemma~\ref{elementary}(i), and $S$ has to contain $R$.
If (2) is satisfied, then $SL_2(R)$ is contained in the finitely generated
group $E_2(S')$ with $S'$ chosen as in Lemm~\ref{elementary}(ii).
\qed

\begin{example} \label{Nagao} \rm
It is a classic result that $SL_n({\mathbb Z}) = E_n({\mathbb Z})$ is
finitely generated for all positive integers $n$. It is also a well-known
result due to Nagao that $SL_2({\mathbb F}_q[t]) = E_2({\mathbb F}_q[t])$
is not finitely generated. However, $SL_2({\mathbb F}_q[t,t^{-1}])$ is
of course finitely generated. More generally, for any $S$--arithmetic
ring (also called a "Hasse domain" in the literature) $R = \cO_S$, the
group $SL_2(\cO_S)$ is finitely generated whenever
the characteristic of $R$ is 0 or the set $S$ of places has cardinality
at least 2. This follows from general results about $S$--arithmetic groups
due to Borel -- Harish--Chandra in characteristic 0 and to Behr in
characteristic $p>0$ (see \cite{BHC} and \cite{Be}).
So Question~\ref{Qu} and, more generally, Question~\ref{RP}
have a positive answer for $S$--arithmetic rings.
\end{example}

In view of the last remark, we are now going to consider rings $R$ with
Krull dimension $> 1$.
For (Laurent) polynomial rings in one variable over ${\mathbb Z}$ or in two
variables over a finite field ${\mathbb F}_q$, Question~\ref{Qu} would involve
the long standing problem mentioned in Remark~\ref{ZLaurent}, which we are not
going to discuss in this paper. So it is natural to consider (Laurent)
polynomial rings in at least two variables over infinite base rings in order
to prove a negative answer to Question~\ref{Qu} for a reasonable class of
rings. Let us fix some further notation:

Let $R_0$ be a finitely generated infinite integral domain with field of
fractions $F_0$. $F$ will be a transcendental extension of $F_0$ of
transcendence degree 2, $F = F_0(s,t)$, and we start by considering
$R = R_0[s,t,s^{-1},t^{-1}]$. (Lemma~\ref{elementary}(ii) indicates that one
should invert the variables in order to avoid trivial obstacles to finite
generation. However, we shall
later see that we can also replace this Laurent polynomial ring with a
polynomial ring.) In view of Lemma~\ref{reformulation}, we want to show that
$SL_2(R)$ is not contained in $E_2(S)$ for any given finitely generated
$R$--subalgebra $S \subseteq F$. The idea is to write $S$ as
${\tilde S}[t^{-1}]$ for a suitable subring ${\tilde S}$ of $S$ and to apply
a similar method as in Section 3 (see Corollary~\ref{notGE2}).
However, it is not enough to just show $SL_2(S) \neq E_2(S)$ here. We need
to be able to exhibit \emph{concrete} elements in $SL_2(S)$ which are not
in $E_2(S)$, and we must be able to choose these elements already in
$SL_2(R)$. So we are going to apply Lemma~\ref{L3} rather than Lemma~\ref{L2}
in the following. The main technical step in the proof is to verify Condition
(1) of Lemma~\ref{L3} in a suitable situation. This will be done in the
framework of the following proposition.

\begin{prop} \label{Mainstep}
Let $f = f(s,t) \in R_0[s,t]$ be a polynomial not divisible by $t$.
Write $f = f_0(s) + f_1(s)t + \ldots f_d(s)t^d$ with polynomials
$f_i \in R_0[s]$ for all $0 \leq i \leq d$. Set
$S := R_0[s,t,s^{-1},t^{-1},f_0^{-1},f^{-1}]$ and $g := 1 - sf_0 \in R_0[s]$.
Let $p$ be any prime element of $R_0$ which does not divide $f_0$ in
$R_0[s]$. Then we obtain
$$\left( \begin{array}{rr} \ 1 - pgt^{-1} & p^2t^{-1} \\
 -g^2t^{-1} & 1 + pgt^{-1} \\ \end{array} \right) \not\in E_2(S)$$
\end{prop}

\proof We first note that $f_0 \neq 0$ since $t$ does not divide $f$. We put
${\tilde S} := R_0[s,t,s^{-1},f_0^{-1}$,$f^{-1}]$ so that $S = {\tilde S}[t^{-1}]$.
Note that $t$ is a prime element of ${\tilde S}$ since it is a prime element of
$R_0[s,t]$ which does not divide $sf_0f$. Replacing $R$ with ${\tilde S}$ and
$\pi$ with $t$, we can now proceed as in Section 3. We introduce the $t$--adic
valuation $v = v_t$ on $F$ with associated discrete valuation ring $\cO$,
maximal ideal $\cP$ and group of units $\cO^*$. We have the same
subgroups $A,B,U$ of $G = SL_2(F)$ as introduced in Fact~\ref{BTtree}. Our
$R[1/\pi]$ is $S$ here, and hence we put $H = SL_2(S)$. Then $A \cap H =
SL_2({\tilde S}), \  B \cap H = D(t^{-1},1)SL_2({\tilde S})D(t,1)$ and
$U \cap H = SL_2\left( \begin{array}{rr} {\tilde S} & {\tilde S}\\
t{\tilde S} & {\tilde S}\\ \end{array} \right)$. We again set
$H_0 = \langle (A \cap H) \cup (B \cap H) \rangle$ and recall that $H_0$
contains $E_2({\tilde S}[t^{-1}]) = E_2(S)$ by Lemma~\ref{1/r}. $g$ and $p$ are
nonzero elements of $R_0[s]$ and hence also elements of $\cO^*$. So we can
consider the matrix $a =  \left(
\begin{array}{rr} 1 & 0\\ g/p & 1 \\ \end{array} \right) \in A = SL_2(\cO)$.

{\bf Claim: $a \not\in HB$.} We introduce another ring, namely
$Z := R_0[s,s^{-1},f_0^{-1}]$. Note that $Z \setminus \{ 0 \} \subseteq \cO^*$
since the elements of $Z$ do not involve $t$. Hence the canonical homorphism
$\phi: {\tilde S} \rightarrow {\tilde S}/t{\tilde S}$ restricted to $Z$ is
injective. $\phi \! \mid_Z$ is also surjective since
$f \equiv f_0 \ {\rm mod} \ t$, implying $\phi(f) = \phi(f_0)$ and
$\phi(f^{-1}) = \phi(f_0^{-1})$. (This was the reason for including $f_0^{-1}$
in $S$.) We can thus decompose the additive group of ${\tilde S}$ as follows:

$(*)\quad {\tilde S} = Z \oplus t{\tilde S} = Z \oplus tZ \oplus t^2{\tilde S}$

We now assume by way of contradiction that $a \in HB$. Then there is a matrix
$h = \left( \begin{array}{rr} \alpha & \beta \\ \gamma & \delta \\ \end{array}
\right) \in H = SL_2(S)$ such that $ha \in B$. Hence we have
$$ \left( \begin{array}{rr} \alpha + \beta g/p & \beta \\
\gamma + \delta g/p & \delta \\ \end{array} \right) \in SL_2
\left( \begin{array}{rr} \cO & t^{-1}\cO \\ t\cO & \cO \\ \end{array} \right)$$

This implies $\delta \in \cO \cap S = {\tilde S}$,
$\beta \in t^{-1}\cO \cap S = t^{-1}(\cO \cap S) = t^{-1}{\tilde S}$,
$\gamma \in \cO \cap S = {\tilde S}$ (recall that $g/p \in \cO^*$) and
$\alpha \in t^{-1}\cO \cap S = t^{-1}{\tilde S}$. Using the decomposition $(*)$,
one therefore finds elements
$a_{-1},a_0,b_{-1},b_0, c_0,c_1,d_0,d_1 \in Z$, $\alpha',\beta' \in t{\tilde S}$
and $\gamma',\delta' \in t^2{\tilde S}$ such that
$\alpha = a_{-1}t^{-1} + a_0 + \alpha'$, $\beta = b_{-1}t^{-1} + b_0 + \beta'$,
$\gamma = c_0 + c_1t + \gamma'$ and $\delta = d_0 + d_1t + \delta'$. But we still
have the conditions $\alpha + \beta g/p \in \cO$ and $\gamma + \delta g/p \in
t\cO$, which yield (together with $p,g \in Z$):
\begin{itemize}
\item[(1)] $pa_{-1} + gb_{-1} = 0$
\item[(2)] $pc_0 + gd_0 = 0$
\end{itemize}
We also have the condition that $\det (h) = \alpha \delta - \beta \gamma = 1$
which leads to the equations $a_{-1}d_0 - b_{-1}c_0 = 0$ (which we do not need)
and
\begin{itemize}
\item[(3)] $a_{-1}d_1 + a_0d_0 - b_{-1}c_1 - b_0c_0 = 1$
\end{itemize}
(which we do need). Since $p$ is prime in $R_0$, it is also prime in $R_0[s]$,
and since $p$ does not divide $f_0$ in $R_0[s]$ by assumption (and certainly not
$s$), $p$ is also a prime element in $R_0[s,s^{-1},f_0^{-1}] = Z$. Furthermore,
$p$ does not divide $g = 1 - sf_0$ in $R_0[s]$, and hence not in $Z$.
Therefore, Equation (1) implies that $p$ divides $b_{-1}$. After cancelling $p$,
the same equation shows that $g$ divides $a_{-1}$. Hence the ideal $(p,g)$ of
$Z$ contains the ideal $(a_{-1},b_{-1})$. Similarly, Equation (2) implies that
$p$ divides $d_0$ in $Z$, then $g$ divides $c_0$ and so $(p,g)$ also contains
the ideal $(c_0,d_0)$ of $Z$. However, (3) implies that
$(a_{-1},b_{-1},c_0,d_0) = (1)$. Hence also $(p,g) = (1)$ in $Z$. This means
that there exist polynomials $x,y \in R_0[s]$ and an integer $n \geq 0$ such
that $px + gy = (sf_0)^n$. Passing to the respective images modulo $p$ and
denoting them by overlining, we obtain ${\overline g}{\overline y} =
({\overline s}{\overline f_0})^n$ in $(R_0/pR_0)[s]$. However, since
$g = 1 - sf_0$, we also have ${\overline g}{\overline z} =
{\overline 1} - ({\overline s}{\overline f_0})^n$ with
$z = 1 + sf_0 + \ldots + (sf_0)^{n-1} \in R_0[s]$. Thus
${\overline g}({\overline y} + {\overline z}) = {\overline 1}$, showing that
${\overline g}$ is a unit in $(R_0/pR_0)[s]$. So ${\overline g} =
{\overline 1} - {\overline s}{\overline f_0} \in (R_0/pR_0)^*$, implying
${\overline s}{\overline f_0} = {\overline 0}$ in $(R_0/pR_0)[s]$. Therefore
$p$ divides $sf_0$ and hence $f_0$ in $R_0[s]$. However, this contradicts
our assumption on $p$. Hence $a \in HB$ is impossible and our claim is proved.

Now we set
$b := \left( \begin{array}{rr}\ 1 & p^2t^{-1} \\ 0 & 1 \\ \end{array} \right)$,
which is certainly an element of $B = SL_2
\left( \begin{array}{rr} \cO & t^{-1}\cO \\ t\cO & \cO \\ \end{array} \right)$
but not of $U = SL_2
\left( \begin{array}{rr} \cO & \cO \\ t\cO & \cO \\ \end{array} \right)$.
Finally we check that $aba^{-1} =
\left( \begin{array}{rr} \ 1 - pgt^{-1} & p^2t^{-1} \\
-g^2t^{-1} & 1 + pgt^{-1} \\ \end{array} \right)$, which is certainly an element
of $H = SL_2(S)$. So by Lemma~\ref{L3}, $aba^{-1} \not\in H_0$, and hence in
particular (since $E_2(S) \leq H_0$ as remarked above)
$aba^{-1} \not\in E_2(S)$.
\qed

\begin{remark} \label{more-examples} \rm
The proof of Proposition~\ref{Mainstep} in fact yields many more matrices
which are not contained in $E_2(S)$. For instance, let $b' =
\left( \begin{array}{rr} \alpha & \beta \\ \gamma & \delta \\ \end{array}
\right)$ be any element of $SL_2(R_0[s])$ satisfying
$\alpha \equiv \delta \ {\rm mod} \ p$, $\beta \equiv 0 \ {\rm mod} \ p^2$
and $\beta \neq 0$. Then $b =
\left( \begin{array}{rr} \alpha & \beta t^{-1} \\ \gamma t & \delta \\
\end{array} \right)$ is an element of $B$, but not of $U$, and again
$aba^{-1} \in SL_2(S)$ (with $a$ as in the proof of
Proposition~\ref{Mainstep}) and $aba^{-1} \not\in E_2(S)$.

One should also note that all these matrices $aba^{-1}$ (including the one
given in Proposition~\ref{Mainstep}) are in fact elements of
$SL_2(R_0[s,t^{-1}])$; their entries only involve $s$ and $t^{-1}$ but not
$s^{-1}$ and $t$.
\end{remark}

We are now in a position to prove the main result of this section, stated as
Theorem~\ref{mainresult} in the Introduction.

\begin{theo} \label{Mainresult}
Let $R_0$ be a finitely generated integral domain with infinitely many
non--associate prime elements, $R = R_0[s,t]$ and $F$ the field of fractions
of $R$. Then $SL_2(R)$ is not contained in $E_2(S)$ for any finitely
generated $R$--subalgebra $S$ of $F$. Equivalently, $\Gamma$ is not a
finitely generated group whenever $SL_2(R) \leq \Gamma \leq SL_2(F)$.
\end{theo}

\proof We first prove the claim for $R' := R_0[s,t^{-1}]$ instead of $R$.
Any finitely generated $R'$--algebra $S' \subseteq F$ is contained in one
of the form $S'' = R_0[s,t,t^{-1},f^{-1}]$ with a nonzero polynomial
$f \in R_0[s,t]$. (Write the generators of $S'$ as
fractions $g_1/f_1,\ldots,g_n/f_n$ with polynomials $g_i,f_i \in
R_0[s,t]$ and define $f$ to be the product $f = f_1\ldots f_n$.) We may
additionally assume that $t$ does not divide $f$ since
$t, t^{-1} \in S''$. Now define $f_0 \in R_0[s] \setminus \{0\}$ as in
Proposition~\ref{Mainstep} and set again
$S = R_0[s,t,s^{-1},t^{-1},f_0^{-1},f^{-1}]$. Since $f_0 \neq 0$, it must
have at least one nonzero coefficient $c$ in $R_0$. Because $c$ is only
divisible by finitely many non--associate prime elements of $R_0$ ($R_0$
is noetherian), our assumption on $R_0$ guarantees the existence of a
prime element $p$ of $R_0$ not dividing $c$ and hence also not dividing
$f_0$ in $R_0[s]$. Now Proposition~\ref{Mainstep} provides us with an
element of $SL_2(R')$ (see the last paragraph of the previous remark)
which is not contained in $E_2(S)$ and hence also not in $E_2(S')$.

So $SL_2(R')$ is not contained in $E_2(S')$ for any finitely generated
$R'$--algebra $S' \subseteq F$. However, the roles of $t$ and $t^{-1}$
are of course symmetric in this situation (more formally: consider the
automorphism of F interchanging $t$ and $t^{-1}$ and fixing $R_0[s]$
pointwise). Hence the analogous statement for $SL_2(R)$ is also true.
Finally, the equivalence of this with the last statement of the theorem
was already established in Lemma~\ref{reformulation}.
\qed

\begin{cor} \label{Maincor}
Let $R$ be one of the following polynomial rings
\begin{itemize}
\item[(1)] $R = \ZZ[t_1,\ldots,t_m]$ with $m \geq 2$
\item[(2)] $R = {\mathbb F}_q[t_1,\ldots,t_m]$ with $m \geq 3$
\end{itemize}
and $F$ its field of fractions. Then any group $\Gamma$ with
$SL_2(R) \leq \Gamma \leq SL_2(F)$ is not finitely generated.

More generally, if $R = R_0'[t_1,\ldots,t_m]$ for an arbitrary finitely
generated integral domain $R_0'$ and $m \geq 3$ or if $R = R_0[t_1,t_2]$
for a Hasse domain (= $S$--arithmetic ring) $R_0$,
then the same result holds for $R$.
(Recall that any Hasse domain has infinitely many non--associate primes
by a well-known theorem from number theory.)
\end{cor}

\begin{remark} \rm
It is interesting to note that for any Hasse domain $R_0$ and any $m > 0$,
$SL_n(R) = E_n(R)$ is finitely generated for $R = R_0[t_1,\ldots,t_m]$
and \emph{all $n \geq 3$}.
This was shown by Suslin in \cite{Sus}.
\end{remark}

I do not know how restrictive the assumption in Theorem~\ref{Mainresult}
concerning the infinitely many
primes really is. It might well be that \emph{any} finitely generated
infinite integral domain has infinitely many non--associate prime elements.
However, I have not yet found a reference yielding this statement in this
generality.

We conclude this section by strengthening the statement $SL_2(S) \neq
E_2(S)$ for finitley generated $R$--subalgebras $S$ of $F$ similarly as
Bachmuth and Mochizuki did for Laurent polynomial rings (see
\cite[Theorem 1]{BM}).

\begin{cor}
Let $R$ and $F$ be as in Theorem~\ref{Mainresult}, and let $S$ be a finitely
generated $R$--subalgebra of $F$. Then any set of generators of $SL_2(S)$
must contain infinitely many elements outside $E_2(S)$.
\end{cor}

\proof If there were a finite subset $L \subset SL_2(S)$ such that
$SL_2(S) = \langle E_2(S) \cup L \rangle$, then there were also a finite
subset $M \subset F$ such that $SL_2(S) \leq E_2(S[M])$ (see the proof of
Lemma~\ref{elementary}(i)), implying $SL_2(R) \leq E_2(S[M])$. Since also
$S[M]$ is a finitely generated $R$--algebra, the latter inclusion is
impossible by Theorem~\ref{Mainresult}.
\qed

\section{Some problems and conjectures} \label{P&C}

One does not necessarily need polynomial rings in at least two variables in
order to get similar results (with, however, more technical proofs) as stated
in Proposition~\ref{Mainstep} and Theorem~\ref{Mainresult}. They all support
the following

\begin{con}\label{Krull3} \rm
If $R$ is a finitely generated integral domain of Krull dimension at least 3
with field of fractions $F$, then a group $\Gamma$ with $SL_2(R) \leq
\Gamma \leq SL_2(F)$ is never finitely generated.
\end{con}

The following question is a natural generalization of Rapinchuk's original
problem in the case $\cG = SL_2$:

\begin{question} \rm
Is it possible in the situation of Conjecture~\ref{Krull3} that there exists
a finitely generated group $\Gamma$ with $SL_2(R) \leq \Gamma \leq SL_2(F')$
if we admit \emph{any} field $F'$ which contains $R$?\\
\end{question}

Whereas the rings of algebraic number theory are very well analyzed
(see Remark~\ref{Nagao}), the situation is pretty unclear for (finitely
generated) domains of \emph {Krull dimension 2}. I think that one only has a
chance to attack Question~\ref{Qu} for this class of rings if one has settled
the following two very concrete (but hard!) problems which were already
mentioned in Remark~\ref{ZLaurent}. I state the first of these two problems as
a conjecture since some numerical evidence (which unfortunately did not lead
to a systematic proof) makes me believe that it is a true statement.

\begin{con} \rm
$SL_2(\ZZ[t,t^{-1}]) = E_2(\ZZ[t,t^{-1}])$.
\end{con}

\begin{question} \rm
Is $SL_2({\mathbb F}_q[t_1,t_1^{-1};t_2,t_2^{-1}])$ generated by
elementary matrices or at least finitely generated?
\end{question}

Rapinchuk's problem for $\cG =SL_n$ with $n \geq 3$ is also challenging but
of a completely different nature. At least for "sufficiently large" $n$ it
purely becomes a question of algebraic $K$-theory, namely the following.

\begin{question} \label{K1} \rm
If $R$ is a finitely generated integral domain, is there always an element
$0 \neq f \in R$ such that $K_1(R[1/f])$ is a finitely generated (abelian)
group?
\end{question}

\begin{remark} \label{n-big} \rm
If Question~\ref{K1} has a positive answer for a given
finitely generated integral domain $R$ with
Krull dimension $d$ and if $n \geq d + 2$, then $SL_n(R[1/f])$ is a
\emph{finitely generated} group containing $SL_n(R)$.

A positive answer to Question~\ref{K1} is known for many rings $R$.
However, it seems to be a hard problem in general. A solution would be
immediately provided by a positive answer to the following more general
question asked by Bass thirty years ago (see \cite{Ba}, problem at the
end of the introduction): Is $K_1(S)$ finitely generated for any
\emph{regular} finitely generated commutative ring $S$? (It is well
known from commutative algebra that for any finitely generated integral
domain $R$, there exists $0 \neq f \in R$ such that $R[1/f]$ is regular.)
\end{remark}

We close this section (and this paper) by returning, in a very special case,
to anisotropic groups, which originally motivated Rapinchuk's
Question~\ref{RP}.

\begin{con} \label{N1} \rm
Let $R$ and $F$ be as in Conjecture~\ref{Krull3}, let $D$ be a quaternion
algebra over $F$, and consider
$\cG = SL_{1,D}$, the (algebraic) group of elements of reduced norm 1.
Then there does not exist a finitely generated group $\Gamma$ with
$\cG(R) \leq \Gamma \leq \cG(F)$.
\end{con}

This conjecture is motivated by Theorem~\ref{Mainresult},
Proposition~\ref{Mainstep} and Remark~\ref{more-examples}. I have no
idea yet what results are to be expected for $SL_{1,D}$ in case the degree
of $D$ is $>2$.


\begin{thebibliography}{99}

\bibitem{BM} {\sc Bachmuth S.}\ and {\sc H. Mochizuki}, $E_2 \neq SL_2$
for most Laurent polynomial rings, {\em American J. of Mathematics}
{\bf 104} (1982), 1181~--~1189.

\bibitem{Ba} {\sc Bass H.}, {\em Introduction to some methods of algebraic
$K$--theory}, Conference Board of the Mathematical Sciences -- Regional
Conference Series in Mathematics {\bf 20}, 1974.

\bibitem{Be} {\sc Behr H.}, Endliche Erzeugbarkeit arithmetischer Gruppen
\"uber Funktionenk\"orpern,
{\em Invent. math.} {\bf 7} (1969), 1~--~32.

\bibitem{BHC} {\sc Borel A.}\ and {\sc Harish--Chandra}, Arithmetic
subgroups of algebraic groups,
{\em Ann. Math.} {\bf 75} (1962), 485~--~535.

\bibitem{Chu} {\sc Chu H.}, On the $GE_2$ of graded rings,
{\em J. Algebra} {\bf 90} (1984), 208~--~216.

\bibitem{Co} {\sc Cohn P.~M.}, On the structure of the $GL_2$ of a ring,
{\em Publ. math. I.H.E.S.} {\bf 30} (1966), 5~--~53.

\bibitem{KS} {\sc Karrass A.} \ and {\sc D. Solitar}, The subgroup of a
free product of two groups with an amalgamated subgroup,
{\em Trans. Am. Math. Soc.} {\bf 150} (1970), 227~--~255.

\bibitem{RSS} {\sc Rapinchuk A.}, \ {\sc Segev Y.} \ and {\sc G. Seitz},
Finite quotients of the multiplicative group of a finite dimensional
division algebra are solvable, {\em J. AMS} {\bf 15} (2002), 929~--~978.

\bibitem{Rap} {\sc Rapinchuk A.}, Algebraic and abstract simple groups:
old and new, Preprint (2002).

\bibitem{Se} {\sc Serre J.-P.}, {\em Trees}, Springer--Verlag, Berlin,
Heidelberg, New York, 1980.

\bibitem{Sus} {\sc Suslin A.} On the structure of the special linear group
over polynomial rings, {\em Isv. Akad. Nauk} {\bf 41} (1977), 235~--~252.

\end{thebibliography}
\end{document}